\documentclass[12pt,a4paper]{article}
\usepackage{amsmath,amssymb,amsthm,bm}
\usepackage{graphicx}
\usepackage{txfonts}

\newtheorem{theorem}{Theorem}[section]
\newtheorem{proposition}[theorem]{Proposition}

\newtheorem{lemma}[theorem]{Lemma}

\newtheorem{example}[theorem]{Example}

\newcommand{\1}{{\bm 1}}

\numberwithin{equation}{section}

\begin{document}

\begin{center}
\large\bf
Complete Multipartite Graphs of Non-QE Class
\end{center}

\bigskip

\begin{center}
Nobuaki Obata\\
Center for Data-driven Science and Artificial Intelligence \\
Tohoku University\\
Sendai 980-8576, Japan \\
obata@tohoku.ac.jp
\end{center}

\bigskip

\begin{quote}
\textbf{Abstract}\enspace
We derive a formula for the QE constant of
a complete multipartite graph
and determine the complete multipartite graphs of non-QE class,
namely, those which do not admit 
quadratic embeddings in Euclidean spaces.
Moreover, we prove that there are exactly four 
primary non-QE graphs among the complete multipartite graphs.
\end{quote}

\begin{quote}
\end{quote}

\begin{quote}
\end{quote}

\section{Introduction}

Realization of a graph $G=(V,E)$ in a Euclidean space $\mathbb{R}^N$
is of fundamental interest.
In this paper, we focus on a particular realization
called quadratic embedding, which traces back to the early works
of Schoenberg \cite{Schoenberg1935,Schoenberg1938} and has been
studied along with Euclidean distance geometry, 
see e.g., \cite{Alfakih2018, Balaji-Bapat2007,
Jaklic-Modic2013, Jaklic-Modic2014, Liberti-Lavor-Maculan-Mucherino2014}.

Let $G=(V,E)$ be a graph, which is always assumed to be
finite and connected. 
A map $\varphi:V\rightarrow \mathbb{R}^N$ is called 
a \textit{quadratic embedding} of $G$ if
\[
\|\varphi(x)-\varphi(y)\|^2
=d_G(x,y),
\qquad
x,y\in V,
\]
where the left-hand side is the square of the Euclidean distance
and the right-hand side is the graph distance.
A graph $G$ is called \textit{of QE class} or
\textit{of non-QE class}
according as it admits a quadratic embedding or not.
It follows from Schoenberg \cite{Schoenberg1935,Schoenberg1938} 
that a graph $G$ admits a quadratic embedding
if and only if the distance matrix $D=[d_G(x,y)]$ is
conditionally negative definite, i.e.,
$\langle f,Df \rangle\le0$ 
for all real column vectors $f$ indexed by $V$
with $\langle \1,f \rangle=0$,
where $\1$ denotes the column vector of which entries are all one
and $\langle\cdot,\cdot\rangle$ is the canonical inner product.
In this connection, a new numeric invariant of a graph
was proposed in the recent papers \cite{Obata2017,Obata-Zakiyyah2018}.
The \textit{quadratic embedding constant} 
(\textit{QE constant} for short) of a graph $G$ is defined by
\begin{equation}\label{1eqn:def od QEC(G)}
\mathrm{QEC}(G)
=\max\{\langle f,Df \rangle\,;\, 
\langle f,f \rangle=1, \, \langle \1,f \rangle=0\},
\end{equation}
where the right-hand side stands for the conditional maximum 
of the quadratic function $\langle f,Df \rangle$
with $f$ running over a unit sphere determined by
$\langle f,f \rangle=1$ and $\langle \1,f \rangle=0$.
By definition, a graph $G$ is of QE class 
if and only if $\mathrm{QEC}(G)\le0$.
An advantage of the QE constant is that
\eqref{1eqn:def od QEC(G)} is determined explicitly
or estimated finely by means of the method of Lagrange's multipliers.
The QE constants are explicitly known 
for some special series of graphs, 
see also \cite{BO-2021,Irawan-Sugeng2021,Lou-Obata-Huang2022, Mlotkowski2022,MO-2018,Obata2023,Purwaningsih-Sugeng2021}.

There are interesting questions both on graphs of QE class 
and on those of non-QE class.
One of the important questions on non-QE graphs 
is to obtain a sufficiently rich list of non-QE graphs.
The main purpose of this paper is to 
determine all complete multipartite graphs of non-QE class.
In this paper we first derive a general formula
for the QE constant of a complete multipartite graph.

\begin{theorem}[Theorem \ref{03thm:main formula}]
\label{01thm:main formula}
Let $k\ge2$ and $m_1\ge m_2\ge \dotsb \ge m_k\ge1$.
The QE constant of the complete $k$-partite 
graph $K_{m_1,m_2,\dots,m_k}$ is given as follows.
\begin{enumerate}
\setlength{\itemsep}{0pt}
\item[\upshape (1)] If $m_1=m_2$, we have
\[
\mathrm{QEC}(K_{m_1,m_2,\dots,m_k})=-2+m_1.
\]
\item[\upshape (2)] If $m_1>m_2$, we have
\[
\mathrm{QEC}(K_{m_1,m_2,\dots,m_k})=-2-\alpha^*,
\]
where $\alpha^*$ is the minimal solution to the equation
\[
\sum_{i=1}^k \frac{m_i}{\alpha+m_i}=0.
\]
Moreover, $-m_1<\alpha^*<-m_2$.
\end{enumerate}
\end{theorem}

Using the above explicit formula,
we determine all complete multipartite graphs with
positive QE constants, that is, of non-QE class.

\begin{theorem}[Theorem \ref{04thm:CMG of non-QE}]
\label{01thm:CMG of non-QE}
Let $k\ge2$ and $m_1\ge m_2\ge \dotsb \ge m_k\ge1$.
Then the complete $k$-partite graph $K_{m_1,m_2,\dots,m_k}$ is 
of non-QE class if and only if one of the following conditions
is satisfied:
\begin{enumerate}
\setlength{\itemsep}{0pt}
\item[\upshape (i)] $m_1\ge3$ and $m_2\ge2$,
\item[\upshape (ii)] $k\ge3$, $m_1\ge5$ and $m_2=\dots=m_k=1$,
\item[\upshape (iii)] $k\ge4$, $m_1=4$ and $m_2=\dots=m_k=1$,
\item[\upshape (iv)] $k\ge5$, $m_1=3$ and $m_2=\dots=m_k=1$.
\end{enumerate}
\end{theorem}

If a graph $H$ is isometrically embedded in a graph $G$,
we have $\mathrm{QEC}(H)\le \mathrm{QEC}(G)$.
Hence, if $G$ contains a non-QE subgraph $H$ isometrically,
then $G$ is of non-QE class too.
Thus, upon classifying graphs of non-QE class
it is important to specify a \textit{primary} non-QE graph,
that is, a non-QE graph $G$ which does not 
contain a non-QE graph $H$ as an
isometrically embedded proper subgraph.

By Theorem \ref{01thm:CMG of non-QE} there are
four families of complete multipartite graphs of non-QE class.
From each family a primary one is specified as follows.

\begin{theorem}[Theorem \ref{04thm:main primary graphs}]
\label{01thm:main primary graphs}
Among the complete multipartite graphs
there are four primary non-QE graphs
which are listed with their QE constants as follows:
\begin{align*}
&\mathrm{QEC}(K_{3,2})=2/5, \\
&\mathrm{QEC}(K_{5,1,1})=1/7, \\
&\mathrm{QEC}(K_{4,1,1,1})=2/7, \\
&\mathrm{QEC}(K_{3,1,1,1,1})=1/7.
\end{align*}
Moreover, any complete multipartite graph of non-QE class
contains at least one of the above four primary non-QE graphs
as an isometrically embedded subgraph.
\end{theorem}

This paper is organized as follows.
In Section 2 we prepare basic notions and notations,
for more details see \cite{MO-2018,Obata-Zakiyyah2018}.
In Section 3 we prove the main formula for the
QE constants of complete multipartite graphs 
(Theorem \ref{01thm:main formula})
and show some examples.
In Section 4 we determine all complete multipartite graphs
of non-QE class (Theorem \ref{01thm:CMG of non-QE})
and specify primary ones (Theorem \ref{01thm:main primary graphs}).
All primary non-QE graphs on six or fewer vertices are
already determined \cite{Obata2023}.
As a result, we find three new primary non-QE graphs on seven vertices.

The QE constant is not only useful for judging whether nor not 
a graph $G$ is of QE-class but also interesting as
a possible scale of classifying graphs \cite{BO-2021}.
Moreover, the QE constant is interesting from the point of
view of spectral analysis of distance matrices,
for distance spectra see e.g., \cite{Aouchiche-Hansen2014, 
Jin-Zhang2014, Lin-Hong-Wang-Shub2013}.
In fact, $\mathrm{QEC}(G)$ lies between
the largest and the second largest eigenvalues of 
the distance matrix $D$
in such a way that $\lambda_2(D)\le \mathrm{QEC}(G)<\lambda_1(D)$.
In this line, an interesting question is to characterize graphs
such that $\lambda_2(D)=\mathrm{QEC}(G)$.
The work is now in progress.

\section{Quadratic Embedding of Graphs}

\subsection{Distance Matrices}

A graph $G=(V,E)$ is a pair of 
a non-empty set $V$ of vertices and a set $E$ of edges,
where $V$ is assumed to be a finite set throughout the paper.
For $x,y\in V$ we write $x\sim y$ if $\{x,y\}\in E$.
A graph is called \textit{connected} if 
any pair of vertices $x,y\in V$ there exists
a finite sequence of vertices $x_0,x_1,\dots,x_m\in V$ 
such that $x=x_0\sim x_1\sim \dotsb\sim x_m=y$.
In that case the sequence of vertices is called
a \textit{walk} from $x$ to $y$ of length $m$.
Unless otherwise stated, by a graph we mean a finite connected graph
throughout this paper.

Let $G=(V,E)$ be a graph.
For $x,y \in V$ with $x\neq y$ let $d_G(x,y)$ denote
the length of a shortest walk connecting $x$ and $y$.
By definition we set $d_G(x,x)=0$.
Then $d_G(x,y)$ becomes a metric on $V$,
which we call the \textit{graph distance}.
The \textit{distance matrix} of $G$ is defined by
\[
D=[d_G(x,y)]_{x,y\in V}\,,
\]
which is a matrix with index set $V\times V$.
For notational simplicity we sometimes 
write $d(x,y)$ for $d_G(x,y)$ when there is no danger of confusion.

Let $C(V)$ be the linear space 
of $\mathbb{R}$-valued functions $f$ on $V$.
We always identify $f\in C(V)$ with a column vector 
$[f_x]_{x\in V}$ through $f_x=f(x)$. 
The canonical inner product on $C(V)$ is defined
\[
\langle f,g \rangle
=\sum_{x\in V} f(x)g(x)\,,
\qquad f,g\in C(V).
\]
The distance matrix $D$ 
induces a linear transformation on $C(V)$
by matrix multiplication as usual.
Since $D$ is symmetric, we have
$\langle f,Dg \rangle=\langle Df,g \rangle$.

\subsection{Quadratic Embedding}

A \textit{quadratic embedding} of a graph $G=(V,E)$ 
in a Euclidean space $\mathbb{R}^N$ is
a map $\varphi:V\rightarrow \mathbb{R}^N$ satisfying
\[
\|\varphi(x)-\varphi(y)\|^2
=d_G(x,y),
\qquad
x,y\in V,
\]
where the left-hand side is the square of the Euclidean distance.
A graph $G$ is called \textit{of QE class} or \textit{of non-QE class}
according as it admits a quadratic embedding or not.

A graph $H=(V^\prime, E^\prime)$ is called a \textit{subgraph}
of $G=(V,E)$ if $V^\prime\subset V$ and $E^\prime\subset E$.
By our convention both $G$ and $H$ are assumed to be connected
and hence admit graph distances of their own.
We say that $H$ is \textit{isometrically embedded} in $G$ if
\[
d_H(x,y)=d_G(x,y),
\qquad x,y\in V^\prime.
\]
In that case we write $H \hookrightarrow G$.

\begin{lemma}\label{02lem:heredity}
Let $G=(V,E)$ and $H=(V^\prime,E^\prime)$ be graphs
and assume that $H$ is isometrically embedded in $G$,
that is, $H \hookrightarrow G$.
\begin{enumerate}
\setlength{\itemsep}{0pt}
\item[\upshape (1)] If $G$ is of QE class, so is $H$.
\item[\upshape (2)] If $H$ is of non-QE class, so is $G$.
\end{enumerate}
\end{lemma}

\begin{proof}
Obvious by definition.
\end{proof}

We say that a graph of non-QE class is \textit{primary}
if it contains no isometrically embedded proper subgraphs
of non-QE class.
In view of Lemma \ref{02lem:heredity},
for classifying graphs of non-QE class
it is essential to explore primary non-QE graphs.

We mention a simple and useful criterion for
isometric embedding.
Recall that the \textit{diameter} of 
a graph $G=(V,E)$ is defined by
\[
\mathrm{diam}(G)=\max\{d_G(x,y)\,;\, x,y\in V\}.
\]

\begin{lemma}\label{02lem:isometric embedding}
Let $G=(V,E)$ be a graph and $H=(V^\prime,E^\prime)$ a subgraph.
\begin{enumerate}
\setlength{\itemsep}{0pt}
\item[\upshape (1)] If $H$ is isometrically embedded in $G$,
then $H$ is an induced subgraph of $G$.
\item[\upshape (2)] If $H$ is an induced subgraph of $G$ and
$\mathrm{diam\,}(H)\le2$, then
$H$ is isometrically embedded in $G$.
\end{enumerate}
\end{lemma}

\begin{proof}
(1) By assumption $d_H(x,y)=d_G(x,y)$ for all $x,y\in V^\prime$.
In particular, for $x,y\in V^\prime$  
$d_H(x,y)=1$ if and only if $d_G(x,y)$=1.
Therefore, if two vertices $x,y\in V^\prime$ are adjacent in $G$,
so are in $H$.
This means that $H$ is an induced subgraph of $G$.

(2) We will prove that
$d_H(x,y)=d_G(x,y)$ for all $x,y\in V^\prime$.
As $d_H(x,y)\le 2$ for $x,y\in V^\prime$ by assumption,
we consider three cases.
If $d_H(x,y)=0$, obviously $d_G(x,y)=0$.
Suppose that $d_H(x,y)=1$.
Then $x$ and $y$ are adjacent in $H$, so are in $G$ and $d_G(x,y)=1$.
Suppose that $d_H(x,y)=2$.
Since $H$ is a subgraph of $G$, 
we have $d_G(x,y)\le d_H(x,y)\le2$.
If $d_G(x,y)=0$, we have $x=y$ and $d_H(x,y)=0$,
which is contradiction.
If $d_G(x,y)=1$, then $x$ and $y$ are adjacent in $G$,
so are in $H$ since $H$ is an induced subgraph of $G$.
We then obtain $d_H(x,y)=1$, which is again contradiction.
Therefore, we have $d_G(x,y)=2$ and $d_H(x,y)=d_G(x,y)$ holds.
\end{proof}

\subsection{Quadratic Embedding Constants}

Let $G=(V,E)$ be a graph with $|V|\ge2$.
The \textit{quadratic embedding constant}
(\textit{QE constant} for short) of $G$ is defined by
\[
\mathrm{QEC}(G)
=\max\{\langle f,Df \rangle\,;\, f\in C(V), \,
\langle f,f\rangle=1, \, \langle \1,f \rangle=0\},
\]
where $\1\in C(V)$ is defined by $\1(x)=1$ for all $x\in V$.

It follows from Schoenberg \cite{Schoenberg1935,Schoenberg1938} 
that a graph $G$ admits a quadratic embedding
if and only if the distance matrix $D=[d_G(x,y)]$ is
conditionally negative definite, i.e.,
$\langle f,Df \rangle\le0$ 
for all $f\in C(V)$ with $\langle \1,f \rangle=0$.
Hence a graph $G$ is of QE class (resp. of non-QE class)
if and only if $\mathrm{QEC}(G)\le0$ (resp. $\mathrm{QEC}(G)>0$).
The idea of QE constant was proposed in 
\cite{Obata-Zakiyyah2018}, where a standard method of
computing the QE constants is established 
on the basis of Lagrange's multipliers.

\begin{lemma}\label{02lem:isometrically embedded subgraphs}
Let $G=(V,E)$ and $H=(V^\prime,E^\prime)$ be two graphs
with $|V|\ge2$ and $|V^\prime|\ge2$.
If $H$ is isometrically embedded in $G$,
we have 
\[
\mathrm{QEC}(H)\le \mathrm{QEC}(G).
\]
\end{lemma}

As the distance matrix of $H$ becomes a principal
submatrix of the distance matrix of $G$,
the proof is straightforward, see also \cite{Obata-Zakiyyah2018}.
Note also that Lemma \ref{02lem:heredity} follows immediately
from Lemma \ref{02lem:isometrically embedded subgraphs}.

For later reference we recall some concrete values of the
QE constants. 
Further examples are found in 
\cite{Irawan-Sugeng2021, Lou-Obata-Huang2022,MO-2018,
Obata2017, Purwaningsih-Sugeng2021}.

\begin{example}[\cite{Obata-Zakiyyah2018}]
\label{02ex:complete graphs}
\normalfont
For the complete graph $K_n$ with $n\ge2$ we have
\[
\mathrm{QEC}(K_n)=-1.
\]
Conversely, any graph $G$ with $\mathrm{QEC}(G)=-1$ is
necessarily a complete graph \cite{BO-2021}.
Moreover, for any graph $G$ on two or more vertices we have
$\mathrm{QEC}(G)\ge-1$.
\end{example}

\begin{example}[\cite{Obata-Zakiyyah2018}]
\label{02ex:cycle graphs}
\normalfont
For the cycles on odd number of vertices we have
\[
\mathrm{QEC}(C_{2n+1})=-\bigg(4\cos^2\dfrac{\pi}{2n+1}\bigg)^{-1},
\qquad n\ge1,
\]
and for those on even number of vertices we have
\[
\mathrm{QEC}(C_{2n})=0,
\qquad
n\ge2.
\]
\end{example}

\begin{example}[\cite{Mlotkowski2022}]
\label{02ex:path graphs}
\normalfont
For the paths $P_n$ with $n\ge2$ we have
\[
\mathrm{QEC}(P_n)=-\bigg(1+\cos\dfrac{\pi}{n}\bigg)^{-1}.
\]
\end{example}

\begin{example}[\cite{Obata-Zakiyyah2018}]
\label{02ex:complete bipartite graphs}
\normalfont
For the complete bipartite graph $K_{m,n}$ with $m\ge1$ and $n\ge1$
we have
\[
\mathrm{QEC}(K_{m,n})=\frac{2(m-1)(n-1)-2}{m+n}\,.
\]
The above formula will be generalized to the complete
multipartite graphs in Subsection 
\ref{03sec:Calculating QE Constants},
see Theorem \ref{03thm:main formula}.
\end{example}

\begin{example}[\cite{Lou-Obata-Huang2022}]
\label{02ex:complete split graph}
\normalfont
For $m\ge1$ and $n\ge1$ the graph join $\bar{K}_m+K_n$ is 
called the \textit{complete split graph}.
We have
\[
\mathrm{QEC}(\bar{K}_m+K_n)=\frac{(m-2)(n-1)-2}{m+n}\,.
\]
For a more general result on the graph join of regular graphs,
see \cite{Lou-Obata-Huang2022}.
\end{example}

\begin{figure}[hbt]
\begin{center}
\includegraphics[width=60pt]{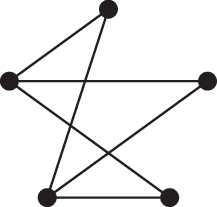}
\quad \raisebox{20pt}{$\cong$} \quad
\includegraphics[width=60pt]{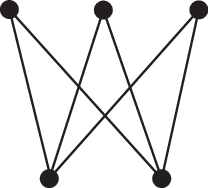}
\qquad
\includegraphics[width=60pt]{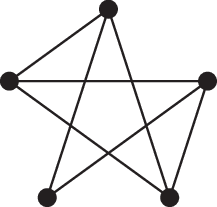}
\quad \raisebox{20pt}{$\cong$} \quad
\includegraphics[width=60pt]{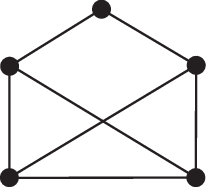}
\end{center}
\caption{Non-QE graphs: G5-10 (left) and G5-17 (right)}
\label{Non-QE graphs on 5 vertices}
\end{figure}

\begin{example}
\label{02ex:small graphs <=5}
\normalfont
A table of the QE constants of graphs on $n\le5$ vertices 
is available \cite{Obata-Zakiyyah2018}.
We see directly from the table that
all graphs on $n\le4$ vertices are of QE-class,
and that there are 21 graphs on five vertices among which
two are of non-QE class.
Those graphs are shown in Figure \ref{Non-QE graphs on 5 vertices},
where G$n$-$x$ stands for the graph on $n$ vertices
with number $x$ in the list of small graphs 
due to McKay \cite{McKay}.
Their QE constants are given by
\[
\mathrm{QEC}(\text{G5-10})=\frac25,
\qquad
\mathrm{QEC}(\text{G5-17})=\frac{4}{11+\sqrt{161}}\approx 0.1689.
\]
Both are primary non-QE graphs since all graphs on four vertices
are of QE class.
Note that G5-10 is the complete bipartite graph $K_{3,2}$.
\end{example}

\section{Complete Multipartite Graphs}

\subsection{Definition and Basic Properties}

For $k\ge2$ let 
$V_1,V_2,\dots V_k$ be mutually disjoint non-empty finite sets.
Setting
\[
V=\bigcup_{i=1}^k V_i,
\qquad
E=\bigcup_{i\neq j}\{\{x,y\}\,;\, x\in V_i,\, y\in V_j\},
\]
we obtain a graph $G=(V,E)$,
which is called a \textit{complete $k$-partite graph} with 
parts $V_1,V_2,\dots,V_k$ and is denoted by $KM(V_1,V_2,\dots,V_k)$.
A \textit{complete multipartite graph} is 
a complete $k$-partite graph for some $k\ge2$.

For $k\ge2$ a complete $k$-partite graph is determined 
(up to graph isomorphisms) by the sequence
$m_1,m_2,\dots,m_k$ defined by $m_i=|V_i|$
and is denoted by
\begin{equation}\label{04def:complete multipartite graph}
K_{m_1,m_2,\dots,m_k}=KM(V_1,V_2,\dots,V_k).
\end{equation}
Without loss of generality we may assume that
$m_1\ge m_2\ge\dotsb\ge m_k\ge1$.
For simplicity we write $K_{1(k)}$ for
$K_{1,1,\dots,1}$ (``$1$'' appears $k$ times).
Obviously, 
\[
K_{1(k)}=K_{1,1,\dots,1}= K_k
\]
which is nothing else but the complete graph on $k$ vertices.
Overusing our notation for the case of $k=1$ we understand
that $KM(V_1)$ is the empty graph on $V_1$, that is the complement
to the complete graph $K_{m_1}$ so that
\begin{equation}\label{03eqn:MV(V_1)}
KM(V_1)= \bar{K}_{m_1}.
\end{equation}
We note that 
\eqref{04def:complete multipartite graph} is not valid for $k=1$.

For later use we list some basic properties
of complete multipartite graphs, 
of which verification is straightforward and is omitted.

\begin{lemma}
Let $k\ge2$ and $m_1\ge m_2\ge \dotsb\ge m_k\ge1$.
Then $\mathrm{diam}(K_{m_1,m_2,\dots,m_k})=2$ 
unless $m_1=m_2=\dots=m_k=1$.
\end{lemma}

\begin{lemma}\label{03lem:graph join structure}
Let $k\ge2$ and let $V_1,V_2,\dots, V_k$ be mutually disjoint
non-empty finite sets.
Consider an arbitrary partition:
\[
\{V_1,V_2,\dots V_k\}
=\{U_1,\dots, U_p\}\cup\{U_1^\prime,\dots,U_q^\prime\},
\]
where $p\ge1$, $q\ge1$ and $p+q=k$.
Then we have
\[
KM(V_1,V_2,\dots V_k)
=KM(U_1,U_2,\dots U_p)+KM(U_1^\prime,U_2^\prime,\dots U_q^\prime),
\]
where the right-hand side is the graph join.
\end{lemma}

For example, we have 
\begin{align*}
K_{m_1,m_2}&=\bar{K}_{m_1}+\bar{K}_{m_2}, \\
K_{m_1,m_2,m_3}&=\bar{K}_{m_2}+K_{m_1,m_3},\\
K_{m_1,m_2,m_3,m_4}&=K_{m_1,m_3}+K_{m_2,m_4},
\quad \text{etc.,}
\end{align*}
where \eqref{03eqn:MV(V_1)} is taken into account.

Next we study a subgraph of $KM(V_1,V_2,\dots V_k)$.
For $2\le p\le k$ let $U_1,U_2,\dots,U_p$ be
$p$ parts chosen from $V_1,V_2,\dots V_k$.
For each $1\le i\le p$ let $U_i^\prime\subset U_i$
be a non-empty subset.
It is then easy to see that
\[
KM(U_1^\prime,U_2^\prime,\dots U_p^\prime)
\hookrightarrow
KM(V_1,V_2,\dots V_k).
\]
This isometric embedding is stated as follows.

\begin{lemma}\label{03lem:isometrically embedded CMG}
Let $m_1\ge m_2\ge \dotsb\ge m_k\ge1$
and $n_1\ge n_2\ge \dotsb\ge n_p\ge1$
be two sequences with $k\ge2$ and $p\ge2$.
If there exist $1\le i_1<i_2<\dots<i_p\le k$ such that
\[
n_1\le m_{i_1},
\quad n_2\le m_{i_2},
\quad \dots,
\quad n_p\le m_{i_p},
\]
then we have
\[
K_{n_1,n_2,\dots,n_p}
\hookrightarrow
K_{m_1,m_2,\dots,m_k}.
\]
\end{lemma}

Obviously, any induced subgraph of
$KM(V_1,V_2,\dots V_k)$ is of the form
$KM(U_1^\prime,U_2^\prime,\dots U_p^\prime)$
as mentioned above.
With the help of Lemma \ref{02lem:isometric embedding}
we come to the following

\begin{lemma}\label{03lem:isometrically embedded subgraphs} 
Let $H$ be a subgraph of a complete
$k$-partite graph $K_{m_1,m_2,\dots,m_k}$.
If $H$ is a graph on two or more vertices and
is isometrically embedded in $K_{m_1,m_2,\dots,m_k}$,
then $H$ is a complete $p$-partite graph 
$K_{n_1,n_2,\dots,n_p}$ of the form as in
Lemma \ref{03lem:isometrically embedded CMG}.
\end{lemma}

\subsection{Calculating QE Constants}
\label{03sec:Calculating QE Constants}

Throughout this subsection,
letting $k\ge2$ and $m_1\ge m_2\ge \dots \ge m_k\ge1$ be
fixed, we set
\[
G=K_{m_1,m_2,\dots,m_k}\,.
\]
Our goal is to obtain an explicit formula for $\mathrm{QEC}(G)$.
Since $1\le \mathrm{diam}(G)\le2$ 
the following general result is useful.

\begin{proposition}{\upshape \cite[Proposition 2.1]{Obata2017}}
\label{03prop:QEC by alpha_min}
Let $G$ be a graph with $1\le \mathrm{diam}(G)\le 2$.
Let $A$ be the adjacency matrix of $G$ and set
\begin{equation}\label{03eqn:def of minimum alpha}
\alpha_{\mathrm{min}}=\min\{\langle f,Af\rangle\,;\,
\langle f,f\rangle=1, \,\,\,
\langle \1,f\rangle=0\}.
\end{equation}
Then we have
\[
\mathrm{QEC}(G)=-2-\alpha_{\mathrm{min}}.
\]
\end{proposition}

Now let $A$ be the adjacency matrix of $G=K_{m_1,m_2,\dots,m_k}$.
Then $A$ is written in a block matrix form:
\begin{equation}\label{03eqn:adjacency matrix in block form}
A=
\begin{bmatrix}
0 & J & J & \dots & J \\
J & 0 & J & \dots & J \\
\vdots &\vdots & \ddots & & \vdots \\
J & J & \dots & 0 & J \\
J & J & \dots & J & 0 
\end{bmatrix}
= J- 
\begin{bmatrix}
J & 0 & 0 & \dots & 0 \\
0 & J & 0 & \dots & 0 \\
\vdots &\vdots & \ddots & & \vdots \\
0 & 0 & \dots & J & 0 \\
0 & 0 & \dots & 0 & J 
\end{bmatrix},
\end{equation}
where $J$ is a matrix of which entries are all one
and the size is understood in the context.
For example, 
in the last block matrix in \eqref{03eqn:adjacency matrix in block form}
$J$'s appear as diagonal entries.
We understand naturally that their sizes are 
$m_1\times m_1,\dots, m_k\times m_k$ from left top to right bottom.

We will find the conditional minimum $\alpha_{\mathrm{min}}$ 
in \eqref{03eqn:def of minimum alpha}.
According to the block matrix form of the adjacency matrix $A$
in \eqref{03eqn:adjacency matrix in block form},
any $f\in C(V)\cong \mathbb{R}^{m_1+\dots+m_k}$ is written as
\[
f=\begin{bmatrix} f_1 \\ f_2 \\ \vdots \\ f_k \end{bmatrix},
\qquad
f_i
=\begin{bmatrix} f_{i1} \\ f_{i2} \\ \vdots \\ f_{im_i}\end{bmatrix}
\in \mathbb{R}^{m_i},
\qquad 1\le i\le k.
\]
Then we have
\begin{equation}\label{03eqn:fAf}
\langle f,Af\rangle
=\langle f, Jf\rangle-\sum_{i=1}^k \langle f_i, Jf_i\rangle
=\langle \1, f\rangle^2-\sum_{i=1}^k \langle \1, f_i\rangle^2.
\end{equation}
Upon applying the method of Lagrange's multiplier we set
\[
\varphi(f,\alpha,\beta)
=\langle f,Af\rangle
 -\alpha(\langle f,f\rangle-1)-\beta\langle \1,f\rangle.
\]
By general theory we know that
the conditional minimum \eqref{03eqn:def of minimum alpha}
is attained at a certain stationary point.
Let $\mathcal{S}$ be the set of all stationary points 
$(f,\alpha,\beta)$, that is, the set of solutions to
\begin{equation}\label{03eqn:def of S}
\frac{\partial\varphi}{\partial f_{ij}}
=\frac{\partial\varphi}{\partial \alpha}
=\frac{\partial\varphi}{\partial \beta}
=0,
\qquad 1\le i\le k,
\quad 1\le j\le m_i\,.
\end{equation}
After simple calculus \eqref{03eqn:def of S} becomes
the following system of equations:
\begin{align}
&(J+\alpha)f_i =-\frac{\beta}{2}\,\1, 
\qquad 1\le i\le k,
\label{04eqn:alpha1} \\
&\langle f,f\rangle=\sum_{i=1}^k \langle f_i,f_i\rangle=1,
\label{04eqn:alpha2} \\
&\langle \1,f\rangle=\sum_{i=1}^k \langle \1,f_i\rangle=0.
\label{04eqn:alpha3}
\end{align}
For $1\le s\le m_i$ the $s$th entry of the left-hand side 
of \eqref{04eqn:alpha1} is given by
\[
((J+\alpha)f_i)_s
=\sum_{t=1}^{m_i}((J)_{st}+\alpha\delta_{st})f_{it}
=\langle \1,f_i\rangle+\alpha f_{is}\,.
\]
Thus \eqref{04eqn:alpha1} is equivalent to
\begin{equation}\label{04eqn:alpha11}
\langle \1,f_i\rangle+\alpha f_{is}=-\frac{\beta}{2}\,,
\qquad 1\le s\le m_i\,,
\quad 1\le i\le k.
\end{equation}
As a result, $\mathcal{S}$ is the set of all solutions
$(f=[f_i],\alpha,\beta)$ to the equations
\eqref{04eqn:alpha11}, \eqref{04eqn:alpha2} and \eqref{04eqn:alpha3}.
For convenience we denote by $\mathcal{S}_A$ the set of
all $\alpha\in\mathbb{R}$ appearing in $\mathcal{S}$,
that is, $(f,\alpha,\beta)\in\mathcal{S}$
for some $f$ and $\beta$.

\begin{lemma}\label{03lem:fAf=alpha}
For any $(f,\alpha,\beta)\in\mathcal{S}$ we have
$\langle f,Af\rangle =\alpha$.
Therefore, 
\[
\alpha_{\mathrm{min}}=\min \mathcal{S}_A\,.
\]
\end{lemma}

\begin{proof}
Let $(f=[f_i],\alpha,\beta)\in \mathcal{S}$.
Noting that the equations 
\eqref{04eqn:alpha1}--\eqref{04eqn:alpha3} are fulfilled,
we have 
\[
\langle f_i\,,\,Jf_i\rangle 
=\bigg\langle f_i\,,\, -\alpha f_i-\frac{\beta}{2}\,\1 \bigg\rangle
=-\alpha \langle f_i\,, f_i\rangle
  -\frac{\beta}{2}\,\langle \1,f_i \rangle
\]
and hence
\[
\sum_{i=1}^k \langle f_i,Jf_i\rangle 
=-\alpha \sum_{i=1}^k \langle f_i, f_i \rangle
  -\frac{\beta}{2}\, \sum_{i=1}^k \langle\1,f_i \rangle
=-\alpha.
\]
On the other hand, by \eqref{04eqn:alpha3} we have
\[
\langle f,Jf\rangle =\langle \1,f\rangle^2=0.
\]
Consequently, we see from \eqref{03eqn:fAf} that
\[
\langle f,Af\rangle
=\langle f, Jf\rangle-\sum_{i=1}^k \langle f_i, Jf_i\rangle
=0-(-\alpha)=\alpha,
\]
as desired.
The assertion in the second half is then apparent.
\end{proof}

In order to study $\mathcal{S}_A$ we come back to the equations
\eqref{04eqn:alpha1}--\eqref{04eqn:alpha3},
where \eqref{04eqn:alpha1} is equivalent to
\eqref{04eqn:alpha11}.

\begin{lemma}\label{03lem:alpha neq0}
If $(\xi_1,\dots,\xi_k,\alpha,\beta)
\in \mathbb{R}^k\times\mathbb{R}\times\mathbb{R}$ satisfies
\begin{align}
(m_i+\alpha)\xi_i &=-\frac{\beta}{2}\,, 
\qquad 1\le i\le k,
\label{03eqn:S1} \\
\sum_{i=1}^k m_i\xi_i^2 &=1, 
\label{03eqn:S2} \\
\sum_{i=1}^k m_i\xi_i &=0,
\label{03eqn:S3}
\end{align}
then, setting $f_i=\xi_i\1$, 
we have $(f=[f_i],\alpha,\beta)\in \mathcal{S}$.
Conversely, 
for any $(f=[f_i],\alpha,\beta)\in \mathcal{S}$ with $\alpha\neq0$
there exist $\xi_1,\dots,\xi_k\in\mathbb{R}$ such that
$f_i=\xi_i\1$ and 
\eqref{03eqn:S1}--\eqref{03eqn:S3} are satisfied.
\end{lemma}

\begin{proof}
It is straightforward to show that 
$(f=[f_i],\alpha,\beta)$ with $f_i=\xi_i\1$
satisfies \eqref{04eqn:alpha11}, \eqref{04eqn:alpha2}
and \eqref{04eqn:alpha3}
under conditions \eqref{03eqn:S1}--\eqref{03eqn:S3}.
We prove the converse.
Suppose that $\alpha\neq0$ and 
$(f=[f_i],\alpha,\beta)$ 
is a solution to
\eqref{04eqn:alpha11}, \eqref{04eqn:alpha2} and \eqref{04eqn:alpha3}.
In view of $\alpha\neq0$ we see from \eqref{04eqn:alpha11} that
$f_{is}$ is constant independently of $1\le s\le m_i$.
In other words, $f_i$ is a constant multiple of $\1$,
say $f_i=\xi_i\1$ with $\xi_i\in\mathbb{R}$.
Then equations \eqref{04eqn:alpha11},
\eqref{04eqn:alpha2} and \eqref{04eqn:alpha3}
are reduced to the equations 
\eqref{03eqn:S1},
\eqref{03eqn:S2} and \eqref{03eqn:S3}, respectively.
\end{proof}

\begin{lemma}\label{03lem:alpha*}
If $\alpha\in\mathcal{S}_A$ and 
$\alpha\not\in\{0,-m_1,-m_2,\dots,-m_k\}$, 
then $\alpha$ satisfies
\begin{equation}\label{3eqn:psi(alpha)}
\sum_{i=1}^k \frac{m_i}{\alpha+m_i}=0.
\end{equation}
Conversely, if $\alpha\in\mathbb{R}$ satisfies \eqref{3eqn:psi(alpha)},
then $\alpha\in\mathcal{S}_A$.
\end{lemma}

\begin{proof}
Let $\alpha\in\mathcal{S}_A$ with
$\alpha\not\in\{0,-m_1,-m_2,\dots,-m_k\}$.
Then $(f,\alpha,\beta)\in\mathcal{S}$ for some $f$ and $\beta$
so that the equations \eqref{04eqn:alpha11},
\eqref{04eqn:alpha2} and \eqref{04eqn:alpha3} are satisfied.
By Lemma \ref{03lem:alpha neq0} 
there exist $\xi_1,\dots,\xi_k\in\mathbb{R}$ such that
$f_i=\xi_i\1$ such that $(f=[f_i],\alpha,\beta)$ fulfills
\eqref{03eqn:S1}--\eqref{03eqn:S3}.
From \eqref{03eqn:S1} we have
\begin{equation}\label{03eqn:in proof Lemma 3.8}
\xi_i=-\frac{\beta}{2}\,\frac{1}{\alpha+m_i}\,,
\qquad
1\le i\le k.
\end{equation}
Then \eqref{03eqn:S2} and \eqref{03eqn:S3} become
\begin{align}
\frac{\beta^2}{4}\sum_{i=1}^k \frac{m_i}{(\alpha+m_i)^2}&=1,
\label{03eqn:in proof Lemma 3.8 (3)} \\
-\frac{\beta}{2}\sum_{i=1}^k \frac{m_i}{\alpha+m_i}&=0,
\label{03eqn:in proof Lemma 3.8 (2)}
\end{align}
respectively.
We see from \eqref{03eqn:in proof Lemma 3.8 (3)}
that $\beta\neq0$ and then from \eqref{03eqn:in proof Lemma 3.8 (2)}
we obtain \eqref{3eqn:psi(alpha)}.

For the converse assertion let $\alpha$ be a solution
to \eqref{3eqn:psi(alpha)}.
Noting that $\alpha\not\in\{-m_1,-m_2,\dots,-m_k\}$,
we define $\xi_i$ by \eqref{03eqn:in proof Lemma 3.8}
and set $f_i=\xi_i\1$.
Choosing $\beta$ satisfying \eqref{03eqn:in proof Lemma 3.8 (3)},
we see that $(f=[f_i],\alpha,\beta)\in \mathcal{S}$
and hence $\alpha\in\mathcal{S}_A$.
\end{proof}

\begin{lemma}\label{03lem:m_1=m_2}
If $m_1=m_2$, then $-m_1\in \mathcal{S}_A$.
\end{lemma}

\begin{proof}
It is straightforward to see that
equations \eqref{03eqn:S1}--\eqref{03eqn:S3} are fulfilled by 
\begin{gather*}
\xi_1=\frac{1}{\sqrt{2m_1}}\,,
\qquad
\xi_2=-\frac{1}{\sqrt{2m_1}}\,,
\qquad
\xi_i=0, \quad 3\le i\le k, 
\\
\alpha=-m_1,
\qquad
\beta=0.
\end{gather*}
It then follows from Lemma \ref{03lem:alpha neq0} that
$(f=[f_i],\alpha,\beta)$ with $f_i=\xi_i\1$,
$\alpha=-m_1$ and $\beta=0$ belongs to $\mathcal{S}$.
Hence $-m_1\in\mathcal{S}_A$.
\end{proof}

\begin{lemma}\label{03lem:m_1>m_2}
If $m_1>m_2$, then $-m_1\not\in\mathcal{S}_A$.
\end{lemma}

\begin{proof}
Suppose that  $-m_1\in\mathcal{S}_A$.
Then $(f,-m_1,\beta)\in\mathcal{S}$ for some $f$ and $\beta$.
Since $\alpha=-m_1\neq0$, by Lemma \ref{03lem:alpha neq0}
there exist $\xi_1,\dots,\xi_k\in\mathbb{R}$
satisfying \eqref{03eqn:S1}--\eqref{03eqn:S3}.
Note that \eqref{03eqn:S1} becomes
\begin{equation}\label{03eqn:S1 in proof(1)}
(m_i-m_1)\xi_i=-\frac{\beta}{2}\,,
\qquad 1\le i\le k.
\end{equation}
Putting $i=1$, we obtain $\beta=0$
and \eqref{03eqn:S1 in proof(1)} becomes
\begin{equation}\label{03eqn:S1 in proof(2)}
(m_i-m_1)\xi_i=0,
\qquad 1\le i\le k.
\end{equation}
Since $m_1>m_2\ge m_3\ge\dotsb\ge m_k\ge1$ by assumption,
we see from \eqref{03eqn:S1 in proof(2)} that $\xi_i=0$ for $2\le i\le k$.
It then follows from \eqref{03eqn:S3} that $\xi_1=0$ too.
Thus we obtain $\xi_i=0$ for all $1\le i\le k$,
which contradicts \eqref{03eqn:S2}.
Consequently, $-m_1\not\in \mathcal{S}_A$.
\end{proof}

We are now in a position to sum up the above argument.

\begin{proposition}\label{03prop:alpha_min determined}
Let $k\ge2$ and $m_1\ge m_2\ge \dotsb \ge m_k\ge1$.
Let $A$ be the adjacency matrix of the complete $k$-partite 
graph $K_{m_1,m_2,\dots,m_k}$ and set
\[
\alpha_{\mathrm{min}}=\min\{\langle f,Af\rangle\,;\,
\langle f,f\rangle=1, \,\,\,
\langle \1,f\rangle=0\}.
\]
\begin{enumerate}
\setlength{\itemsep}{0pt}
\item[\upshape (1)] If $m_1=m_2$, then $\alpha_{\mathrm{min}}=-m_1$.
\item[\upshape (2)] If $m_1>m_2$, then $\alpha_{\mathrm{min}}=\alpha^*$,
where $\alpha^*$ is the minimal solution to 
\begin{equation}\label{3eqn:psi(alpha)(1)}
\sum_{i=1}^k \frac{m_i}{\alpha+m_i}=0.
\end{equation}
Moreover, $-m_1<\alpha^*<-m_2$.
\end{enumerate}
\end{proposition}

\begin{proof}
We set
\[
\psi(\alpha)=\sum_{i=1}^k \frac{m_i}{\alpha+m_i}\,.
\]
It follows from Lemma \ref{03lem:alpha*} that
\begin{equation}\label{03eqn:inclusion}
\mathcal{S}_A
\subset \{0,-m_1,\dots,-m_k\}
\cup\{\alpha\in\mathbb{R}\,;\,\psi(\alpha)=0\},
\end{equation}
where $\{\alpha\in\mathbb{R}\,;\,\psi(\alpha)=0\}=\emptyset$
may occur.
Since $\psi(\alpha)<0$ for all $\alpha<-m_1$,
we have $\{\alpha\in\mathbb{R}\,;\,\psi(\alpha)=0\}
\subset (-m_1,+\infty)$.
Therefore, the minimum of the right-hand side 
of \eqref{03eqn:inclusion} is $-m_1$.

(1) By Lemma \ref{03lem:m_1=m_2} we know that
$-m_1\in \mathcal{S}_A$.
We then see from \eqref{03eqn:inclusion} that
$\min \mathcal{S}_A=-m_1$ so that $\alpha_{\mathrm{min}}=-m_1$.

(2) We first note from $m_1>m_2$ that $\psi(\alpha)=0$ has a solution.
Hence the minimal solution $\alpha^*$ certainly exists
and,  as is easily verified by simple calculus, we have
$-m_1<\alpha^*<-m_2$.
Moreover, $\alpha^*\in\mathcal{S}_A$ by Lemma 3.8.
On the other hand, $-m_1\not\in \mathcal{S}_A$ by 
Lemma \ref{03lem:m_1>m_2}.
Therefore,
$\min \mathcal{S}_A=\alpha^*$ so that $\alpha_{\mathrm{min}}=\alpha^*$.
\end{proof}

We now come to the first main result.

\begin{theorem}[Theorem \ref{01thm:main formula}]
\label{03thm:main formula}
Let $k\ge2$ and $m_1\ge m_2\ge \dotsb \ge m_k\ge1$.
\begin{enumerate}
\setlength{\itemsep}{0pt}
\item[\upshape (1)] If $m_1=m_2$, 
we have $\mathrm{QEC}(K_{m_1,m_2,\dots,m_k})=-2+m_1$.
\item[\upshape (2)] If $m_1>m_2$, 
we have $\mathrm{QEC}(K_{m_1,m_2,\dots,m_k})=-2-\alpha^*$,
where $\alpha^*$ is the minimal solution to \eqref{3eqn:psi(alpha)(1)}.
Moreover, $-m_1<\alpha^*<-m_2$.
\end{enumerate}
\end{theorem}

\begin{proof}
Straightforward from Propositions \ref{03prop:QEC by alpha_min}
and \ref{03prop:alpha_min determined}.
\end{proof}

\subsection{Some Special Cases}

We apply Theorem \ref{03thm:main formula} to some special cases.

\begin{example}[See Example \ref{02ex:complete bipartite graphs}]
\normalfont
The QE constant of a bipartite graph $K_{m,n}$ is
obtained by simple algebra.
If $m=n\ge1$, we have 
\begin{equation}\label{03eqn:QEC(K_mm)}
\mathrm{QEC}(K_{m,m})=-2+m.
\end{equation}
If $m>n$, we have
\begin{equation}\label{03eqn:QEC(K_mn)}
\mathrm{QEC}(K_{m,n})=\frac{2(m-1)(n-1)-2}{m+n}\,.
\end{equation}
The former formula \eqref{03eqn:QEC(K_mm)} is reproduced
by setting $m=n$ in the latter \eqref{03eqn:QEC(K_mn)}.
Moreover, since the right-hand side of \eqref{03eqn:QEC(K_mn)}
is symmetric in $m$ and $n$, we see that
\eqref{03eqn:QEC(K_mn)} is valid for any $m\ge1$ and $n\ge1$.
\end{example}

\begin{example}
\normalfont
The QE constant of a tripartite graph $K_{l,m,n}$ is
obtained by simple algebra.
If $l=m\ge n\ge1$, we have 
\begin{equation}\label{03eqn:QEC(K_lln)}
\mathrm{QEC}(K_{l,l,n})=-2+l.
\end{equation}
Consider the case of $l>m\ge n\ge1$.
The equation \eqref{3eqn:psi(alpha)(1)} becomes
\[
\frac{l}{\alpha+l}+\frac{m}{\alpha+m}+\frac{n}{\alpha+n}
=0
\]
and the minimal solution is given by
\[
\alpha^*
=-\frac{lm+mn+nl+\sqrt{(lm+mn+nl)^2-3lmn(l+m+n)}}{l+m+n}\,.
\]
Then we have 
\begin{equation}\label{03eqn:QEC(K_lmn)}
\mathrm{QEC}(K_{l,m,n})
=-2+\frac{lm+mn+nl+\sqrt{(lm+mn+nl)^2-3lmn(l+m+n)}}{l+m+n}\,.
\end{equation}
By simple algebra we see that
\eqref{03eqn:QEC(K_lln)} is reproduced from
\eqref{03eqn:QEC(K_lmn)} by setting $l=m$.
Moreover, since the right-hand side of \eqref{03eqn:QEC(K_lmn)} 
is symmetric in $l,m$ and $n$,
we conclude that the formula \eqref{03eqn:QEC(K_lmn)}
is valid for any $l\ge1$, $m\ge1$ and $n\ge1$.
\end{example}

\begin{example}[see Example \ref{02ex:complete split graph}]
\label{03ex:complete split graph}
\normalfont
Let $m\ge1$ and $n\ge1$.
The $(n+1)$-partite graph $K_{m,1(n)}$ coincides with the graph join
$\bar{K}_m+K_n$, which is known as a \textit{complete split graph}.
The QE constant is obtained by simple algebra.
If $m=1$, we have 
\begin{equation}\label{03eqn:QEC(K_11(n))}
\mathrm{QEC}(K_{1,1(n)})=-2+1=-1.
\end{equation}
In fact, $K_{1,1(n)}=K_{n+1}$ is the complete graph
and the result is mentioned also in Example \ref{02ex:complete graphs}.
Assume that $m>1$.
The minimal solution to 
\[
\frac{m}{\alpha+m}+\frac{n}{\alpha+1}
=0
\]
is given by $\alpha^*=-m(n+1)/(m+n)$ and we then come to
\begin{equation}\label{03eqn:QEC(K_m1(n)}
\mathrm{QEC}(K_{m,1(n)})=-2-\alpha^*
=-2+\frac{m(n+1)}{m+n}=\frac{(m-2)(n-1)-2}{m+n}\,.
\end{equation}
Note that \eqref{03eqn:QEC(K_11(n))} is reproduced from
\eqref{03eqn:QEC(K_m1(n)} by setting $m=1$.
Hence the formula \eqref{03eqn:QEC(K_m1(n)} 
is valid for any $m\ge1$ and $n\ge1$.
\end{example}

\begin{example}
\normalfont
For $n\ge2$ the complete multipartite graph
$K_{2(n)}$ is called a \textit{cocktail party graph}
or a \textit{hyperoctahedral graph}.
A straightforward application of Theorem \ref{03thm:main formula}
yields
\[
\mathrm{QEC}(K_{2(n)})=0.
\]
Note that $K_{2(n)}$ is obtained by deleting $n$ disjoint
edges from the complete graph $K_{2n}$.
It is known \cite{Obata-Zakiyyah2018} that
the QE constant of 
a graph obtained by deleting two or more disjoint edges
from a complete graph is zero.
\end{example}

\section{Primary Non-QE Graphs}

\subsection{Complete Multipartite Graphs of Non-QE Class}

With the help of Theorem \ref{03thm:main formula} 
(Theorem \ref{01thm:main formula} in Introduction)
all complete multipartite graphs of non-QE class are determined.

Let $k\ge2$ and $m_1\ge m_2\ge \dotsb \ge m_k\ge1$.
If $m_1\ge3$ and $m_2\ge2$, 
then $K_{m_1,m_2,\dots,m_k}$ is of non-QE class
because $K_{3,2}\hookrightarrow K_{m_1,m_2,\dots,m_k}$
and $K_{3,2}$ is of non-QE class.
We will examine the rest cases.

(Case 1) $m_1\ge3$ and $m_2=1$.
In that case we have $m_1>m_2=\dots=m_k=1$
so that $K_{m_1,m_2,\dots,m_k}=K_{m_1,1(k-1)}$ becomes a
complete split graph.
Using the formula \eqref{03eqn:QEC(K_m1(n)}
in Example \ref{03ex:complete split graph}, we obtain
\[
\mathrm{QEC}(K_{m_1,1(k-1)})=\frac{(m_1-2)(k-2)-2}{m_1+k-1}\,.
\]
Hence $\mathrm{QEC}(K_{m_1,1(k-1)})>0$ if and only if
one of the following three conditions are fulfilled:
(i) $k\ge 3$ and $m_1\ge5$;
(ii) $k\ge 4$ and $m_1\ge4$;
(iii) $k\ge 5$ and $m_1\ge3$.

(Case 2) $m_1=2$.
If $m_1=m_2$, we see from Theorem \ref{03thm:main formula} that
\[
\mathrm{QEC}(K_{2,2,m_3,\dots,m_k})=-2+m_1=0.
\]
If $m_1>m_2$, we have $m_2=\dots=m_k=1$.
Applying the formula \eqref{03eqn:QEC(K_m1(n)}, we have
\[
\mathrm{QEC}(K_{2,1(k-1)})=-\frac{2}{k+1}<0.
\]

(Case 3) $m_1=1$. In that case we have $m_2=\dots=m_k=1$
so that $K_{1(k)}=K_k$ becomes a complete graph.
We know that $\mathrm{QEC}(K_{1(k)})=-1$.

Thus, there are no non-QE graphs in Cases 2 and 3. 
Summing up the above argument, we claim the following

\begin{theorem}[Theorem \ref{01thm:CMG of non-QE}]
\label{04thm:CMG of non-QE}
Let $m_1\ge m_2\ge \dotsb \ge m_k\ge1$ with $k\ge2$.
Then the complete $k$-partite graph $K_{m_1,m_2,\dots,m_k}$ is 
of non-QE class if and only if one of the following conditions
is satisfied:
\begin{enumerate}
\setlength{\itemsep}{0pt}
\item[\upshape (i)] $m_1\ge3$ and $m_2\ge2$,
\item[\upshape (ii)] $k\ge3$, $m_1\ge5$ and $m_2=\dots=m_k=1$,
\item[\upshape (iii)] $k\ge4$, $m_1=4$ and $m_2=\dots=m_k=1$,
\item[\upshape (iv)] $k\ge5$, $m_1=3$ and $m_2=\dots=m_k=1$.
\end{enumerate}
\end{theorem}

The four families of complete multipartite graphs
in Theorem \ref{04thm:CMG of non-QE} are mutually exclusive.
With the help of Lemma \ref{03lem:isometrically embedded subgraphs},
from each family we find a primary non-QE graph.
The result is stated in the following 

\begin{theorem}[Theorem \ref{01thm:main primary graphs}]
\label{04thm:main primary graphs}
Among the complete multipartite graphs
there are four primary non-QE graphs,
which are listed with their QE constants as follows:
\begin{align*}
&\mathrm{QEC}(K_{3,2})=2/5, \\
&\mathrm{QEC}(K_{5,1,1})=1/7, \\
&\mathrm{QEC}(K_{4,1,1,1})=2/7, \\
&\mathrm{QEC}(K_{3,1,1,1,1})=1/7.
\end{align*}
Moreover, any complete multipartite graph of non-QE class
contains at least one of the above four primary non-QE graphs
as an isometrically embedded subgraph.
\end{theorem}

\subsection{Exploring Primary Non-QE graphs}

In Theorem \ref{04thm:main primary graphs} we find
three primary non-QE graphs on seven vertices.
Note also that they are all complete split graphs:
\begin{align*}
K_{5,1,1}&=\bar{K}_5+K_2=\text{G7-351},\\
K_{4,1,1,1}&=\bar{K}_4+K_3=\text{G7-774},\\
K_{3,1,1,1,1}&=\bar{K}_3+K_4=\text{G7-845},
\end{align*}
where G$n$-$x$ stands for the graph on $n$ vertices
with number $x$ in the list of small graphs due to McKay \cite{McKay}.

\begin{figure}[hbt]
\begin{center}
\includegraphics[width=60pt]{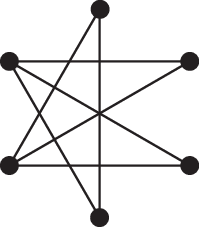}
\qquad\qquad
\includegraphics[width=60pt]{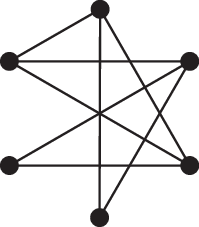}
\qquad\qquad
\includegraphics[width=60pt]{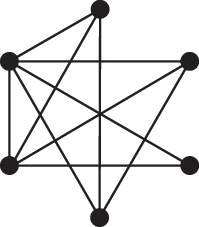}
\end{center}
\caption{Primary non-QE graphs:
G6-30, G6-60 and G6-84 (from left to right)}
\label{04fig:Primary non-QE graphs on 6 vertices}
\end{figure}

On the other hand, 
as is mentioned in Example \ref{02ex:small graphs <=5},
there are exactly two primary non-QE graphs on five vertices.
Moreover, it is proved \cite{Obata2023} that 
there are exactly three primary non-QE graphs on six vertices.
They are the graphs G6-30, G6-60 and G6-84,
see Figure \ref{04fig:Primary non-QE graphs on 6 vertices}.
For the readers' convenience, we record their QE constants:
\begin{align*}
\mathrm{QEC}(\text{G6-30})
&=\frac{-4+\sqrt{19}}{3}\approx 0.1196, \\
\mathrm{QEC}(\text{G6-60})
&=\lambda_{60}^*\approx 0.2034, \\
\mathrm{QEC}(\text{G6-84})
&=\lambda_{84}^*\approx 0.1313,
\end{align*}
where $\lambda_{60}^*$ is
the unique positive root of
$5\lambda^3+26\lambda^2+24\lambda-6=0$,
and $\lambda_{84}^*$ is the unique positive root of
$3\lambda^4+14\lambda^3+18\lambda^2+5\lambda-1=0$.

In this line an interesting question is to systematically 
construct primary non-QE graphs
and to study the distribution of their QE constants.

\bigskip
{\bfseries Acknowledgements:}
This work is supported by
JSPS Grant-in-Aid for Scientific Research No.~19H01789.
The author is grateful to the referees for their useful comments 
which improved the paper.


\end{document}